\documentclass[12pt,oneside,english]{amsart}
\usepackage[latin1]{inputenc}
\usepackage{amssymb}

\makeatletter
\newtheorem{theorem}{Theorem}

\newtheorem{example}{Example}

\newtheorem{remark}{Remark}

\usepackage{babel}

\makeatother
\begin{document}
\baselineskip=17pt

\title[A recursion for divisor function]{A recursion for divisor function over divisors belonging to a prescribed finite sequence of positive integers and a solution of the Lahiri problem for divisor function $\sigma_x(n)$}

\author{Vladimir Shevelev}
\address{Departments of Mathematics \\Ben-Gurion University of the
 Negev\\Beer-Sheva 84105, Israel. e-mail:shevelev@bgu.ac.il}


\begin{abstract}
For a finite sequence of positive integers $A=\{a_j\}_{j=1}^{k},$ we prove a recursion for divisor function $\sigma_{x}^{(A)}(n)=\sum_{d|n,\enskip d\in A}d^x.$  As a corollary, we give an affirmative solution of the problem posed in 1969 by D. B. Lahiri [3]: to find an identity for divisor function $\sigma_x(n)$ similar to the classic
pentagonal identity in case of $x=1.$
 \end{abstract}

\maketitle

\section{Introduction and main results}
    We start with the two well known beautiful classical recursions. Let $p(n)$ be the number of all partitions of positive integer $n$ and $\sigma(n)$ be the sum of its divisors. Then (sf [1],[5]) we have
 \begin{equation}\label{1}
 p(n)=p(n-1)+p(n-2)-p(n-5)-p(n-7)+p(n-12)+p(n-15)-...
  \end{equation}
  \begin{equation}\label{2}
 \sigma(n)=\sigma(n-1)+\sigma(n-2)-\sigma(n-5)-\sigma(n-7)+\sigma(n-12)+\sigma(n-15)-...
  \end{equation}
  where the numbers 1,2,5,7,12,15,... appearing in the successive terms in (1)-(2) are the positive pentagonal numbers $\{v_m\}$
  given by
  \begin{equation}\label{3}
  v_m=m(3m\mp1)/2,\enskip m=1,2,...
  \end{equation}
  In identities (1)-(2) we accept that $p(m)=0,\enskip\sigma(m)=0$ when $m<0.$ The only formal difference is that
  (1) is true with the understanding that
   \begin{equation}\label{4}
  p(0)=1,
  \end{equation}
   while (2) is valid with the understanding that
    \begin{equation}\label{5}
    \sigma(0)=n.
    \end{equation}
  Note that, formulas (1)-(2) are proved with help of the famous Euler pentagonal identity
  \begin{equation}\label{6}
  \prod_{n=1}^{\infty}(1-q^n)=\sum_{m=-\infty}^{\infty}(-1)^{m}q^{m(3m-1)/2}.
  \end{equation}
  In its turn, a combinatorial proof of (6) is based on the following statement (sf [1]). Let $p_e(n)\enskip (p_o(n))$ denote the number of partitions of $n$ into even (odd) number of distinct parts. Then
  \begin{equation}\label{7}
  p_e(n)-p_o(n)=\begin{cases} (-1)^m,\;\;if\;\;n=m(3m\mp1)/2,\\0, \;\; otherwise\end{cases}.
   \end{equation}
  Let $\sigma_{x}(n)$ denote the sum of the $x$th powers of the divisors of $n.$
  In 1969, Lahiri [3] noticed that every definition of $\sigma_k(0)=f(n), \enskip k\neq1$ is irrelevant in order to keep the classical identity (2) and posed the following problem: "Whether analogous identities exist for divisor function  $\sigma_k(n)$ of higher degree?" Formally, for every not necessarily integer value of $x,\enskip -\infty<x<\infty,$  for $\sigma_{x}(n)$  we could consider an identity of type (2)
  of the form
  $$ \sigma(n)=g_x(n)+\sigma(n-1)+\sigma(n-2)-\sigma(n-5)-\sigma(n-7)+...,$$
  where $ \{g_x(n)\}$ is some "compensating sequence," and a solution of the Lahiri problem consists of a description
  of the compensating sequence for every $n$ without a reference to its divisors.
   In particular, by the definition of $\sigma_{x}(n),$ and accepting as in (5) $ \sigma(0)=n,$ we find
  $$ g_x(1)=0, \enskip g_x(2)=2^x-2,\enskip g_x(3)=3^x-2^x-1, \enskip g_x(4)=4^x-3^x-1,$$
   $$g_x(5)=5^x-4^x-3^x-2^x+4, \enskip g_x(6)=6^x-5^x-4^x,... $$
  At first sight, this sequence is even more complicated than $\sigma_{x}(n),$ and it seems hardly probable to find a required description of it. Our paper, in particular, is devoted to this aim. For a simplification of our transformations, below we accept the unique convention
   \begin{equation}\label{8}
\sigma_x(n)=0, \enskip if \enskip n\leq0
 \end{equation}
 It is easy to see that in this case we have only a little change of the compensating sequence in the identity of the same form
$$\sigma_x(n)=h_x(n)+\sigma_x(n-1)+\sigma_x(n-2)-\sigma_x(n-5)-\sigma_x(n-7)+... ,$$
such that
$$ h_x(n)=g_x(n)+\begin{cases} (-1)^{m-1},\;\;if\;\;n=m(3m\mp1)/2,\\0, \;\; otherwise\end{cases}.$$ 
Note that this relation is so simple only due to Euler pentagonal identity (6); in more general case (see below Theorem 1) the corresponding relations could be very complicated and the convention (8) plays the unique role for the obtaining of general result. In particular, we write (1)-(2) in just a little another form. Namely, according to (8), instead of conventions (4)-(5), we accept the unique convention
$$p(0)=0, \enskip \sigma(0)=0.$$\newpage
Then with help of (7) it is easy to see that, instead of (1)-(2), we have
 $$ p(n)=h^{(p)}(n)+$$
 \begin{equation}\label{9}
 p(n-1)+p(n-2)-p(n-5)-p(n-7)+p(n-12)+p(n-15)-...,
  \end{equation}
  where the compensating sequence has the form
 \begin{equation}\label{10}
   h^{(p)}(n)=\begin{cases} (-1)^{m-1},\;\;if\;\;n=m(3m\mp1)/2,\\0, \;\; otherwise\end{cases}.
    \end{equation}
    and, in view of the same structure of (1) and (2) and taking into account (4)-(5), we see that
  $$\sigma(n)=h^{(\sigma)}(n)+$$
  \begin{equation}\label{11}
 \sigma(n-1)+\sigma(n-2)-\sigma(n-5)-\sigma(n-7)+\sigma(n-12)+\sigma(n-15)-...,
  \end{equation}
  where
 \begin{equation}\label{12}
   h^{(\sigma)}(n)=\begin{cases} (-1)^{m-1}n,\;\;if\;\;n=m(3m\mp1)/2,\\0, \;\; otherwise\end{cases}.
    \end{equation}

    Before formulating a generalization of (9) and (11), we study the divisor function over divisors belonging to a prescribed finite sequence $A$ of positive integers. In the trivial case of a one-element sequence $A=\{a\}$ we
    put
    \begin{equation}\label{13}
\sigma_{x}^{(\{a\})}(n)= \begin{cases} a^x,\;\;if\;\;a|n,\enskip n>0,\\0, \;\; otherwise\end{cases},\enskip x\in (-\infty, +\infty).
 \end{equation}
According to (13), we  accept
\begin{equation}\label{14}
\sigma_{x}^{(\{a\})}(n)=0,\enskip n\leq0,
\end{equation}
such that
\begin{equation}\label{15}
\sigma_{x}^{(\{a\})}(n)=\sigma_{x}^{(\{a\})}(n-a)+\begin{cases} a^x,\;\;if\;\;n=a,\\0, \;\; otherwise\end{cases}.
\end{equation}
Consider now, for a fixed $k\geq1,$ an arbitrary sequence
 \begin{equation}\label{16}
 A=\{a_j\}_{j=1}^k
  \end{equation}
 of positive integers. For a fixed $x,$ let us consider an associated sequence
\begin{equation}\label{17}
B(A; x)=\{b_i(x)\}_{i=1}^{2^k},
 \end{equation}
where
\begin{equation}\label{18}
b_i(x)=a_{j_1}^x+a_{j_2}^x+a_{j_3}^x+...+a_{j_r}^x,
\end{equation}
if the binary expansion of $i-1$ is
\begin{equation}\label{19}
i-1=2^{j_1-1}+2^{j_2-1}+...+2^{j_r-1},\enskip 1\leq j_1<j_2<...<j_r,\enskip 1\leq r\leq k.
\end{equation}
In particular, since $2^k-1=2^{1-1}+2^{2-1}+...+2^{k-1},$ then\newpage
\begin{equation}\label{20}
b_{2^k}(x)=a_{1}^x+a_{2}^x+...+a_{k}^x,
\end{equation}
while, since to $i=1$ corresponds the empty set of terms in (19), then
\begin{equation}\label{21}
b_1(x)=0.
\end{equation}
Furthermore,
\begin{equation}\label{22}
b_2(x)=a_1^x,\enskip b_3(x)=a_2^x,\enskip b_4(x)=a_1^x+a_2^x,\enskip etc.
\end{equation}
Moreover, denote
\begin{equation}\label{23}
b_i(1)=b_i,\enskip 1\leq i\leq 2^k.
\end{equation}
For $n\geq1,$ consider divisor function over sequence $A$
\begin{equation}\label{24}
\sigma_{x}^{(A)}(n)=\sum_{d|n,\enskip d\in A}d^x
\end{equation}
in the understanding that every term $d^x$ repeats correspondingly to the multiplicity of $d$ in sequence $A.$
Besides, we accept the convention
\begin{equation}\label{25}
\sigma_{x}^{(A)}(n)=0, \enskip if \enskip n\leq0.
\end{equation}
Denote by $ \{t_n \}$ the  Thue-Morse sequence [4], [2] which is defined as
\begin{equation}\label{26}
t_n=(-1)^{s(n)},
\end{equation}
where $s(n)$ denotes the number of ones in the binary expansion of $n.$
\begin{theorem}\label{1}
In convention  $\sigma (n\leq0)=0, $ we have the following recursion
\begin{equation}\label{27}
\sigma_{x}^{(A)}(n)=h^{(A)}_x(n)+\sum_{i=2}^{2^k}t_{2i-1}\sigma_{x}^{(A)}(n-b_i)
\end{equation}
where the compensating sequence $h^{(A)}_x(n)$ is defined as
\begin{equation}\label{28}
h^{(A)}_x(n)=\sum_{i\geq2:\enskip b_i=n}t_{2i-1}b_i(x).
\end{equation}
\end{theorem}
\begin{remark}
Taking into account that
$$1+s(i-1)=s(2(i-1)+1)=s(2i-1),$$
we prefer to write $t_{2i-1}$ instead of $\enskip-t_{i-1}.$
\end{remark}
Note that, as follows from (28), for $n>b_{2^k}, \enskip h^{(A)}_x(n)=0$ such that
\begin{equation}\label{29}
\sigma_{x}^{(A)}(n)=\sum_{i=2}^{2^k}t_{2i-1}\sigma_{x}^{(A)}(n-b_i),\enskip n> b_{2^k}.
\end{equation}
Consider now the divisor function\newpage
\begin{equation}\label{30}
\sigma_{x}(n)=\sum_{d|n}d^x.
\end{equation}
Putting here
\begin{equation}\label{31}
b_i(x)=j_1^x+j_2^x+j_3^x+...+j_r^x \enskip(and \enskip b_i=b_i(1))
\end{equation}
if the binary expansion of $i-1$ is defined by (19),
we obtain the following result.
\begin{theorem}\label{2}
We have
$$\sigma_{x}(n)=h_x(n)+$$
\begin{equation}\label{32}
\sigma_x(n-1)+\sigma_x(n-2)-\sigma_x(n-5)-\sigma_x(n-7)+\sigma_x(n-12)+\sigma_x(n-15)-...,
\end{equation}
where the  compensating sequence $\{h_x(n)\}$ is defined as
\begin{equation}\label{33}
h_x(n)=\sum_{i\geq2:\enskip b_i=n}t_{2i-1}b_i(x),\enskip n\geq1.
\end{equation}
\end{theorem}
Theorem 2 gives a solution of the  Lahiri problem  for divisor function $\sigma_x(n).$

\section{Proof of Theorem 1}
We use the induction over the number of elements of sequence $A,$ the base of which is given by (15). Note that if, instead of $A=\{a_1,...,a_k\},$ to consider
the sequence
\begin{equation}\label{34}
A'=\{a_1,...,a_k, a_{k+1}\},
\end{equation}
  then we have
\begin{equation}\label{35}
\sigma_{x}^{(A')}(n)=\sigma_{x}^{(A)}(n)+\sigma_{x}^{(\{a_{k+1}\})}(n).
\end{equation}
Furthermore, in the case of $A',$ to every $i,\enskip 1\leq i\leq 2^k,$ with the binary expansion (19) of $i-1$ corresponds bijectively the number $2^k+i$ from $[2^k+1,2^{k+1}]$ with the expansion
$$ 2^k+i-1=2^{j_1-1}+2^{j_2-1}+...+2^{j_r-1}+2^k$$
such that the associated sequence has the form
\begin{equation}\label{36}
b_i(x)=\begin{cases}a_{j_1}^x+a_{j_2}^x+a_{j_3}^x+...+a_{j_r}^x,\;\;if\;\;1\leq i\leq 2^k,\\a_{j_1}^x+a_{j_2}^x+a_{j_3}^x+...+a_{j_r}^x+a_{k+1}^x, \;\;if\;\;2^k+1\leq i\leq 2^{k+1}.\end{cases}
\end{equation}
This means that, for $1\leq l\leq 2^k,$ we have
\begin{equation}\label{37}
b_{l+2^k}(x)=b_l(x)+a_{k+1}^x\enskip (in \enskip particular,\enskip b_{1+2^k}(x)=a_{k+1}^x).
\end{equation}
Notice also, that\newpage
\begin{equation}\label{38}
t_{2^{k+1}+1}=1; \enskip t_{2(l+2^k)-1}=t_{2l+2^{k+1}-1}=-t_{2l-1}
\end{equation}
and
\begin{equation}\label{39}
\sum_{1\leq l\leq 2^k}t_{2l-1}=-\sum_{1\leq l\leq 2^k}t_{l-1}=0.
\end{equation}
Suppose now that the theorem is true up to $k.$ Then, using (37)-(38), we have
$$\sum_{i=2}^{2^{k+1}}t_{2i-1}\sigma_{x}^{(A')}(n-b_i)=$$
$$\sum_{i=2}^{2^{k}}t_{2i-1}\sigma_{x}^{(A')}(n-b_i)+\sum_{i=2^k+1}^{2^{k+1}}t_{2i-1}\sigma_{x}^{(A')}(n-b_i)=$$
$$\sum_{i=2}^{2^{k}}t_{2i-1}\sigma_{x}^{(A')}(n-b_i)+
\sum_{l=1}^{2^{k}}t_{2(l+2^k)-1}\sigma_{x}^{(A')}(n-b_{l+2^k})=$$
\begin{equation}\label{40}
\sum_{i=2}^{2^{k}}t_{2i-1}\sigma_{x}^{(A')}(n-b_i)-
\sum_{l=1}^{2^{k}}t_{2l-1}\sigma_{x}^{(A')}((n-b_l)-a_{k+1}).
\end{equation}
Furthermore, by (40) and (35), we have
$$\sum_{i=2}^{2^{k+1}}t_{2i-1}\sigma_{x}^{(A')}(n-b_i)=$$ $$\sum_{i=2}^{2^{k}}t_{2i-1}\sigma_{x}^{(A)}(n-b_i)+
\sum_{i=2}^{2^{k}}t_{2i-1}\sigma_{x}^{(\{a_{k+1}\})}(n-b_i)$$
\begin{equation}\label{41}
-\sum_{i=1}^{2^{k}}t_{2i-1}\sigma_{x}^{(A)}((n-b_i)-a_{k+1})
-\sum_{i=1}^{2^{k}}t_{2i-1}\sigma_{x}^{(\{a_{k+1}\})}((n-b_i)-a_{k+1}).
\end{equation}
Note that, according to (15),
$$\sum_{i=2}^{2^{k}}t_{2i-1}\sigma_{x}^{(\{a_{k+1}\})}(n-b_i)-
\sum_{i=1}^{2^{k}}t_{2i-1}\sigma_{x}^{(\{a_{k+1}\})}((n-b_i)-a_{k+1})$$
\begin{equation}\label{42}
=\sigma_{x}^{(\{a_{k+1}\})}(n)+a_{k+1}^x\sum_{1\leq i\leq 2^{k}:\enskip n-b_i=a_{k+1}}t_{2i-1}.
\end{equation}
Therefore, from (41) we find
$$\sum_{i=2}^{2^{k+1}}t_{2i-1}\sigma_{x}^{(A')}(n-b_i)=\sigma_{x}^{(\{a_{k+1}\})}(n)+a_{k+1}^x\sum_{1\leq i\leq 2^{k}:\enskip n-b_i=a_{k+1}}t_{2i-1}+$$\newpage
\begin{equation}\label{43}
\sum_{i=2}^{2^{k}}t_{2i-1}\sigma_{x}^{(A)}(n-b_i)-\sum_{i=1}^{2^{k}}t_{2i-1}\sigma_{x}^{(A)}((n-a_{k+1})-b_i),
\end{equation}
or, using the inductive hypothesis, we have
$$\sum_{i=2}^{2^{k+1}}t_{2i-1}\sigma_{x}^{(A')}(n-b_i)=\sigma_{x}^{(\{a_{k+1}\})}(n)+a_{k+1}^x\sum_{1\leq i\leq 2^{k}:\enskip n-b_i=a_{k+1}}t_{2i-1}-$$
\begin{equation}\label{44}
(\sigma_{x}^{(A)}((n-a_{k+1})-h^{(A)}(n-a_{k+1}))+\sigma_x^{(A)}(n)-h^{(A)}(n).
\end{equation}
Furthermore,
$$\sum_{2\leq i\leq 2^{k+1}:\enskip b_i=n}t_{2i-1}b_i(x)=\sum_{2\leq i\leq 2^{k}:\enskip b_i=n}t_{2i-1}b_i(x)+$$
$$\sum_{2^k+1\leq i\leq 2^{k+1}:\enskip b_i=n}t_{2i-1}b_i(x)=h_x^{(A)}-\sum_{1\leq l\leq 2^{k}:\enskip b_{l+2^k=n}}t_{2l-1}b_{l+2^k}(x)=$$
\begin{equation}\label{45}
h_x^{(A)}(n)-h_x^{(A)}(n-a_{k+1})-a_{k+1}^x\sum_{1\leq l\leq 2^{k}:\enskip b_{l}=n-a_{k+1}}t_{2l-1}.
\end{equation}
Finally, summing the results of (44) and (45), we complete our proof:
$$\sum_{i=2}^{2^{k+1}}t_{2i-1}\sigma_{x}^{(A')}(n-b_i)+\sum_{2\leq i\leq 2^{k+1}:\enskip b_i=n}t_{2i-1}b_i(x)=$$
$$\sigma_{x}^{(\{a_{k+1}\})}(n)+\sigma_x^{(A)}(n)=\sigma_x^{(A')}(n).\blacksquare$$

\section{Proof of Theorem 2}
If to consider as a finite sequence $A$ the sequence $A=A_k=\{1,2,...,k\},$ then, for $n\leq k,$ we have
\begin{equation}\label{46}
\sigma_x^{(A_k)}(n)=\sigma_x(n)
\end{equation}
and, by Theorem 1, the $(\pm)$-structure of $\sigma_x^{(A_k)}(n)$ is the same as in the case of $x=1$ (see (11)). Therefore, independently from the summands (either $\sigma_1(n)$ or $\sigma_x(n)$) we have the same reductions, i.e.
$$\sigma_x(n)=\sigma_x^{(A_k)}(n)=h_x^{(A_k)}(n)+$$
$$\sigma_x^{(A_k)}(n-1)+\sigma_x^{(A_k)}(n-2)-\sigma_x^{(A_k)}(n-5)-\sigma_x^{(A_k)}(n-7)+$$
\begin{equation}\label{47}
\sigma_x^{(A_k)}(n-12)+\sigma_x^{(A_k)}(n-15)-...,\enskip (n\leq k),
\end{equation}
with the compensating sequence
\begin{equation}\label{48}
h_x^{(A_k)}(n)=\sum_{i\geq2:\enskip b_i=n}t_{2i-1}b_i(x)
\end{equation}\newpage
where  $b_i(x)$ are defined by (31).
If, instead of $A_k,$ to consider $N,$ then for every $n$ we actually consider a finite part of (47) which corresponds to
$A_n=\{1,2,...,n\}.$ Thus (47) is true for $A=N,$ and (32)-(33) follow.$\blacksquare$ \newline
\begin{example}\label{1}
Consider the case of $x=1,$ i.e. the case of sum-of-divisors function.\end{example}
 Then we have
$$h_1^{(N)}(n)=n\sum_{i\geq2:\enskip b_i=n}t_{2i-1}=-n\sum_{i\geq2:\enskip b_i=n}(-1)^{s(i-1)}=n(p_o(n)-p_e(n))$$
and, in view of (7), we obtain (11) as a special case of Theorem 2.
\section{Expression of compensating sequence $\{h_x(n)\}=\{h_x(n)^{(N)}\}$ via known sequences }
Note that from the definition of sequence $b_n(x)$ (see (31) and (19)) it follows that if
\begin{equation}\label{49}
n-1=\sum_{i\geq1}\beta(i)2^{i-1}
\end{equation}
is the binary expansion of $n-1,$ then
\begin{equation}\label{50}
b_n(x)=\sum_{i\geq1}\beta(i)i^x,
\end{equation}
such that
\begin{equation}\label{51}
b_n=\sum_{i\geq1}\beta(i)i.
\end{equation}
Notice that, (51) is Sequence A029931(n-1) in [6]). Denoting
\begin{equation}\label{52}
A029931(n)=\eta(n),
\end{equation}
 according to (33), we have
\begin{equation}\label{53}
h_x(n)=\sum_{j\geq1:\enskip \eta(j)=n}(-1)^{s(j)-1}b_{j+1}(x).
\end{equation}
\begin{example}\label{2}
Consider the case of $x=0,$ i.e. the case of the number of divisors of $n.$
\end{example}
Then, by (50) and (53), the compensating sequence has the form
\begin{equation}\label{54}
h_0(n)=\sum_{j\geq1:\enskip \eta(j)=n}(-1)^{s(j)-1}s(j),
\end{equation}
where $s(n),$ as in the above, is the number of ones in the binary expansion of $n.$ The first terms of compensating sequence $\{h_0^{(N)}(n)\}_{n\geq1}$ are:\newpage
\begin{equation}\label{55}
1,1, -1,-1,-3,0,-2,1,2,1,2,4,1,-1,...
\end{equation}
Note that, in view of (49), the expression (51) gives the number of all partitions with distinct parts of a fixed
values of $b_n.$ This means that if to denote by $e_n$  the set of the terms of A029931 for which $\eta(j)=n:$
$$ e_1=\{1\}, e_2=\{2\}, e_3=\{3,4\}, e_4=\{5,8\}, e_5=\{6,9,16\},$$ $$e_6=\{7,10,17,32\}, e_7=\{11,12,18.33.64\} \enskip ..., $$
then the concatenation of this sets leads to the ordering of all partitions of n with distinct parts :
$$n=\sum_{i\geq1}\beta(i)i$$
 respectively to the values of $\sum_{i\geq1}\beta(i)2^{i-1}.$ Thus this way leads us to the Adams-Watters sequence "Decimal equivalent of binary encoding of partitions into distinct parts" (see A118462 in [6]):
\begin{equation}\label{56}
1,2,3,4,5,8,6,9,16,7,10,17,32,11,12,18,33,64,...
\end{equation}
Denote Sequence (56) via $W(n)$ and put $|e_n|=R(n),$ where, for $n\geq1,\enskip \{R(n)\}$  is Sequence A000009[6] that is the number of partitions of $n$ into distinct parts. Finally, denote
\begin{equation}\label{57}
T(n)=\sum_{k=1}^{n}R(k)=A036469(n)-1,\enskip n\geq1.
\end{equation}
Then from (53) we find
\begin{equation}\label{58}
h_x(n)=\sum_{m=0}^{R(n)-1}(-1)^{s(W(T(n)-m))-1}b_{W(T(n)-m)+1}(x).
\end{equation}
\begin{example}\label{3}
Let us calculate the seventh term $h_0(7)$ of sequence (55).
\end{example}
By (58), here we have
\begin{equation}\label{59}
h_0(n)=\sum_{m=0}^{R(n)-1}(-1)^{s(W(T(n)-m))-1}s(W(T(n)-m)).
\end{equation}
If $n=7,$ then we have from the corresponding tables of [6]:
$$R(n)=5,\enskip T(n)=18, \enskip W(18)=64,\enskip  W(17)=33, $$
$$\enskip W(16)=18,\enskip  W(15)=12, \enskip W(14)=11.$$
Thus, according to (59), we find
$$h_0(7)=1-2-2-2+3=-2.$$\newpage
\section{Some another identities }
   In case of the finite set $$A_k=\{2^{j-1}\}_{j=1}^k,$$ according to (18), we have
   \begin{equation}\label{60}
b_i(x)= 2^{(j_1-1)x}+2^{(j_2-1)x}+...+2^{(j_r-1)x},
\end{equation}
if
\begin{equation}\label{61}
i-1=2^{j_1-1}+2^{j_2-1}+...+2^{j_r-1},\enskip 1\leq j_1<j_2<...<j_r,\enskip 1\leq r\leq k.
\end{equation}
Thus $b_i$ has an especially simple form:
\begin{equation}\label{62}
b_i=i-1.
\end{equation}
Let
$$n=2^{\alpha_1-1}+...+2^{\alpha_m-1}.$$
Then, according to Theorem 1 and (62), we find
$$h_{x}^{(A_{k})}(n)=
\sum_{i\geq 2 :\enskip b_i=n}t_{2i-1}b_i(x)=$$
\begin{equation}\label{63}
t_{2n+1}b_{n+1}(x)=(-1)^{s(2n+1)}(2^{(\alpha_1-1)x}+...+2^{(\alpha_m-1)x})
\end{equation}
and
$$\sigma_{x}^{(A_{k})}(n)=(-1)^{s(2n+1)}(2^{(\alpha_1-1)x}+...+2^{(\alpha_m-1)x})+$$
\begin{equation}\label{64}
\sum_{i\geq 2}(-1)^{s(2i-1)}\sigma_{x}^{(A_k)}(n-(i-1)).
\end{equation}
Considering now the infinite sequence of powers of 2:
$$A=\{2^{j-1}\}_{j\geq1},$$
we conclude that
 $$\sigma_{x}^{(A)}(n)=(-1)^{s(2n+1)}(2^{(\alpha_1-1)x}+...+2^{(\alpha_m-1)x})+$$
$$\sum_{i\geq 2}^{n}(-1)^{s(2i-1)}\sigma_{x}^{(A)}(n-(i-1)),$$
or
$$\sigma_{x}^{(A)}(n)=(-1)^{s(2n+1)}(2^{(\alpha_1-1)x}+...+2^{(\alpha_m-1)x})-$$
\begin{equation}\label{65}
\sum_{j=1}^{n-1}(-1)^{s(n-j)}\sigma_{x}^{(A)}(j)),
\end{equation}
where
$$n=2^{\alpha_1-1}+...+2^{\alpha_m-1}.$$
In particular, in the case of $x=0,$ we obtain the identity
\begin{equation}\label{66}
\sum_{j=1}^{n}(-1)^{s(n-j)}\sigma_{0}^{(A)}(j))=(-1)^{s(n)-1}s(n).
\end{equation}
The sequence $\{\sigma_{0}^{(A)}(n))-1\}_{n\geq1}$ is well-known so-called "the binary carry sequence" (A007814 in [6]). In the case of $x=1,$ we obtain the identity
\begin{equation}\label{67}
\sum_{j=1}^{n}(-1)^{s(n-j)}\sigma_{1}^{(A)}(j))=(-1)^{s(n)-1}n.
\end{equation}
The sequence $\{\sigma_{1}^{(A)}(n))\}_{n\geq1}$ is also well-known (see A038712 in[6]).


\begin{thebibliography}{6}
\bibitem 1 G.\enskip E.\enskip Andrews,  \enskip \slshape The theory of partitions, \upshape Addison-Wesley, 1976.
\bibitem 2 S.\enskip Goldstein, K.\enskip A.\enskip Kelly and E.\enskip R.\enskip Speer,\enskip \slshape The fractal structure of rarefied sums of the Thue-Morse sequence, \upshape J. of Number Theory \bfseries 42 \mdseries (1992), 1-19.
\bibitem 3 D.\enskip B. \enskip Lahiri,\enskip \slshape Identities connecting elementary divisor function of different degrees, and allied congruences,\upshape\enskip Math. Scand., \enskip\bfseries 24\mdseries\enskip(1969),\enskip102--110.
\bibitem 4 M.\enskip Morse,\enskip \slshape Reccurent geodesics on a surface of negative
curvature,\upshape \enskip Trans. Amer. Math. Soc.\bfseries \enskip 22 \mdseries(1921),\enskip  84-100.
\bibitem 5  I. \enskip Niven and H.\enskip S.\enskip Zuckerman \slshape An introduction to the theory of numbers,\upshape\enskip John Wiley, New York, 1960.
\bibitem 6 N.\enskip J.\enskip A.\enskip Sloane,\enskip\slshape The On-Line Encyclopedia of Integer Sequences \upshape(http: //www.research.att.com)

\end{thebibliography}
\end{document}